\documentstyle[12pt]{article}

  \newcommand{\const}{\rm const}

  \newcommand{\Dom}{\rm  Dom}

\begin{document}

 \begin{center}

 \ {\bf Relations between  growth of entire functions}\par

\vspace{3mm}

 \ {\bf and behavior of its Taylor coefficients.}

\vspace{5mm}

{\bf  M.R.Formica, E.Ostrovsky, and L.Sirota. } \par

 \vspace{5mm}

 \end{center}

 \ Universit\`{a} degli Studi di Napoli Parthenope, via Generale Parisi 13, Palazzo Pacanowsky, 80132,
Napoli, Italy. \\

e-mail: mara.formica@uniparthenope.it \\

\vspace{3mm}

Department of Mathematics and Statistics, Bar-Ilan University,
59200, Ramat Gan, Israel. \\

\vspace{4mm}

e-mail: eugostrovsky@list.ru\\
Department of Mathematics and Statistics, Bar-Ilan University,\\
59200, Ramat Gan, Israel.

\vspace{4mm}

e-mail: sirota3@bezeqint.net \\

\vspace{5mm}

\begin{center}

 \ {\it Abstract.}

 \end{center}

 \ We derive in the closed and unimprovable form the bilateral non - asymptotic relations between  growth of entire functions
 and  decay rate at infinity of its Taylor coefficients. We investigate the functions of one as well as of several complex variables. \par
 \ We will apply the convex analysis: Young - Fenchel (Legendre) transform, Young inequality, saddle - point method etc. \par

 \vspace{5mm}

 \par

\vspace{5mm}

\section{Notations, definitions, statement of problem, previous results.}

\vspace{5mm}

 \ Let $ \  f = f(z) \ $ be entire (analytical) complex valued function defined on the whole complex plane:

 \begin{equation} \label{source function}
 f(z) = \sum_{k=0}^{\infty} c_k \ z^k, \ \hspace{4mm} c_k = c_k[f],
 \end{equation}
i.e. such that the radius of convergence of the  Taylor (power) series (\ref{source function}) is equal to infinity:

$$
\overline{\lim}_{n \to \infty} \sqrt[n]{|c_n|} = 0.
$$

 \ Recall that the so - called  {\it maximal function,} or equally {\it maximal majorant} \ $ \ M_f(r) = M(r), \ r \in (0,\infty) \ $
for the source one is defined as follows

\begin{equation} \label{maximal fun}
M_f(r) \stackrel{def}{=} \max_{z: |z| \le r} |f(z)| = \max_{z: |z| = r} |f(z)|.
\end{equation}

\vspace{4mm}

 \ {\it It is the one of the classical problem from the theory of the entire functions: establish the relations between the
 asymptotical behavior as $ \ n \to \infty \ $ the coefficients $ \ c_n = c_n[f] \ $ and the asymptotical behavior for the maximal
 function $ \ M_f(r) = M(r) \ $ as} $ \ r \to \infty. \ $ \par

 \vspace{3mm}

 \ See for example the classical monographs \cite{Evgrafov},  \cite{Levin}; where are described also some  important applications, for instance
 in the theory of distributions of zeros of entire functions and in the functional analysis, in particular in the theory of operators.\par
 \ As regards for the modern work we mention an article \cite{Freund}. \par

\vspace{4mm}

 \ Let us bring some notions from the classical theory.  The {\it order }  $ \ \rho[f] \ $ of the function $ \ f = f(z) \ $
may be calculated by the  following relations

\vspace{4mm}

\begin{equation} \label{rho f}
\rho[f] = \overline{\lim}_{n \to \infty} \left\{ \ \frac{n \ \ln n}{|\ln \ | c_n \ | \ |} \  \right\}.
\end{equation}

\vspace{4mm}

 \ Correspondingly, the {\it type} $ \ \beta[f] \ $ is equal to

\vspace{4mm}

\begin{equation} \label{beta f}
\beta[f] = \overline{\lim}_{n \to \infty} \left\{ \ n^{1/\rho[f]} \ \sqrt[n]{ |c_n|} \  \right\}.
\end{equation}

\vspace{4mm}

\ {\bf Our target in this short report is to establish the non - asymptotical exact bilateral estimations between
the coefficients of the considered function and its maximal majorant. \par
 \ They have a very simple, closed and very general form. \par
 \ We extend obtained results in the last section onto the entire functions of several complex variables.} \par

 \vspace{3mm}

{\it \ We derive as a consequence the natural conditions for the possible coincidence of these estimates: upper and lower ones.} \par

 \vspace{4mm}

  \ Recall also one important for us auxiliary facts  from the theory of convex function.
  The so - called Young - Fenchel, or Legendre transform $ \ g^*(y), \ y \in R \ $  for the
 given numerical valued function $ \ g = g(x) \ $ having the  convex non - empty domain of definition $ \ \Dom[g] \ $ is defined as follows

\begin{equation} \label{Young Fen}
g^*(y) \stackrel{def}{=} \sup_{x \in \Dom[g]} (xy - g(x)).
\end{equation}

 \ The function $ \  g^*(\cdot)  \ $ named as customary {\it Young conjugate,} or simple {\it conjugate} to the source one $ \ g(\cdot). \ $ \par

 \ Many  examples and properties of this transform with applications e.g. to the theory of Orlicz spaces may be found in the classical monographs
\cite{KrasRut}, \cite{RaoRen1}, \cite{RaoRen2}, \cite{Rockafellar}, \cite{Tiel}.

\vspace{3mm}

 \ The famous theorem of Fenchel - Moreau tell us that if the function $ \ g = g(x) \ $ is convex and continuous defined on the convex set,
 then

\begin{equation} \label{Fen Mor}
g^{**}(x) = g(x),
\end{equation}

see \cite{Rockafellar}, chapters 2,3; \cite{Tiel}.  \par

 \ Recall also the following important Young's inequality

\begin{equation} \label{Young in}
xy \le g(x) + g^*(y), \ x,y \in R;
\end{equation}

or more generally

\begin{equation} \label{Young in gamma}
xy \le g(\gamma x) + g^*(y/\gamma), \  \gamma = \const > 0, \ x,y \in R.
\end{equation}

\vspace{5mm}

\section{Main results: upper and lower non - asymptotic estimates.}

\vspace{5mm}

\begin{center}

 \ {\sc Upper estimate.} \par

\end{center}

\vspace{4mm}

\ We mean to derive the non - asymptotic {\it upper} estimate of the Taylors coefficients  for the function
$ \ f: \ c_n[f] \ $ via its maximal function $ \ M_f(r). \ $ \par

 \vspace{4mm}

 \ Introduce the function

\begin{equation} \label{Lamfun}
\Lambda(v) \stackrel{def}{=} \ln M_f(e^v), \ v \in R.
\end{equation}

\vspace{3mm}

 \ However, it is sufficient for us hereinafter to suppose at last instead (\ref{Lamfun}) only an
 unilateral restriction

\begin{equation} \label{unilateral}
M_f(r) \le \exp \left( \  \Lambda(\ln r)  \ \right), \ r > 0,
\end{equation}
for certain function $ \ \Lambda  =\Lambda(v), \ \in R. \ $ The restrictions on this function will be clarified below. \par

\vspace{4mm}

 \ {\bf Proposition 2.1.}   We propose under condition (\ref{unilateral})

\begin{equation} \label{upp cn}
|c_n| \le \exp \left( \  - \Lambda^*(n)  \  \right), \ n = 1,2,\ldots.
\end{equation}

\vspace{4mm}

 \ {\bf Proof.} We start from the simple  estimate

$$
|c_n| \le \frac{M_f(r)}{r^n}, \ r > 0.
$$

 \ Therefore

 $$
 |c_n| \le \exp \left( \  - n \ln r + \ln M_f(r)  \ \right)  = \exp \left( \ -nv + \Lambda(v)   \ \right) =
 $$

$$
\exp \left( \ -(nv - \Lambda(v))   \ \right),  \ v = \ln r \in (-\infty, \infty).
$$

 \ Since the last inequality is true for arbitrary value $ \ v, \ $ one can take the minimum over $ \ v: \ $

$$
|c_n| \le \exp ( \  - \sup_{v \in R} (nv - \Lambda(v))  \  ) = \exp \left( \ - \Lambda^*(n)   \ \right),
$$
Q.E.D. \par

\vspace{4mm}

\begin{center}

 \ {\sc Lower estimate.} \par

\end{center}

\vspace{4mm}

 \ We intent  to derive in this subsection  the non - asymptotic {\it lower} estimate
 of the Taylors coefficients $ \ c_n[f] = c_n \ $ for the function $ \ f: \ c_n[f] \ $ via its maximal
 function $ \ M_f(r). \ $ \par
 \ Equivalently: we want to obtain the {\it upper} estimate for the maximal function $ \ M_f(r) \ $
 through its series of Taylor coefficients.\par

\vspace{4mm}

 \ The lower estimate is more complicated. We will follow  the authors on an article \cite{Kozachenko}. Suppose

\begin{equation} \label{lower cn}
|c_n| \le \exp \left( \ - Q(n)  \  \right)
\end{equation}
for  certain increasing to infinity function $ \ Q = Q(z), \ z \ge 1. \ $   We have for all the sufficient greatest values $ \ r \ge e \ $

\begin{equation} \label{Mr est1}
M_f(r) \le  \sum_{n=0}^{\infty} \exp \left( \ n v - Q(n)  \ \right) = : R(v) \ R_Q(v), \ v = \ln r.
\end{equation}

 \ It follows from the theory of the saddle - point method, see e.g. \cite{Fedoryuk}, chapters 1,2, that for some finite
 positive constant $ \ C = C[Q] \in (1,\infty) \ $

$$
R_Q(v) \le \exp( \ \sup_{n}(  Cnv - Q(n)) \ ) =
$$

\begin{equation} \label{saddle}
 \exp( \ Q^{*}(C v) \ ) = \exp( \ Q^{*}(C \ \ln r) \ ), \ r \ge e.
\end{equation}

 \ If in addition  the function $ \ Q = Q(v) \ $ is in turn {\it Young conjugate} to some  continuous and convex one,
 say $ \  G(\cdot): \ Q = G^*,  \ $ then by virtue of Theorem of Fenchel - Moreau

\begin{equation} \label{prelim Mr est}
M_f(r) \le   C_0[G] \exp \left( \  G(C \ \ln r)  \ \right),  \  r \ge e.
\end{equation}

 \ In particular, the function $ \ G(\cdot)  \ $  may coincides with the introduced before function $ \ \Lambda^*(\cdot); \ $
 and we conclude in this case

\begin{equation} \label{prelim Mr lambda}
M_f(r) \le C_0[\Lambda] \exp \left( \  \Lambda(C \ \ln r)  \ \right),  \  r \ge e.
\end{equation}

 \vspace{4mm}

  \ Let us bring the strong proof of (\ref{prelim Mr est}), of course, under appropriate conditions. Define the following function

$$
K_Q(\epsilon) := \sum_{n=0}^{\infty} \exp \left( \ -\epsilon Q^*(n) \ \right), \ \epsilon \in (0,1).
$$
 \ It is proved in  particular in  \cite{Kozachenko} that if $ \ \exists \epsilon \in (0,1) \ \Rightarrow K(\epsilon) < \infty,  \  $ then

$$
R(v) \le K_Q(\epsilon) \exp \left( \ (1 - \epsilon) Q^{**} \left(  \  \frac{v}{1 - \epsilon} \ \right)    \ \right).
$$
 \ Further, as long as

$$
n v \le Q^*( ( 1 - \epsilon) \ v) + Q^{**}(n/(1 - \epsilon)), \ \epsilon \in (0,1).
$$
we have

$$
R(v) \le U_Q(\epsilon) \exp \left( \ Q^{**}(v/(1 - \epsilon)) \ \right),
$$
where in addition to the (\ref{saddle})

$$
U(\epsilon) =  U_Q(\epsilon) := \sum_{n=0}^{\infty} \exp (Q^*((1 - \epsilon) \ n) - Q^*(n)).
$$

\vspace{3mm}

\ To summarize:

\vspace{4mm}

 \ {\bf Proposition 2.2.} Denote

$$
Y_Q(\epsilon) = Y(\epsilon) := \min(U(\epsilon), K(\epsilon)), \ \epsilon \in (0,1).
$$

 \ We assert

\begin{equation} \label{General lower}
R_Q(v) \le Y_Q(\epsilon) \ \exp \left( \ Q^{**}(v/(1 - \epsilon)) \ \right),  \ \epsilon \in (0,1).
\end{equation}

 \ Of course,

\begin{equation} \label{General lower  inf}
R_Q(v) \le  \inf_{\epsilon \in (0,1)} \left\{ \ Y_Q(\epsilon) \ \exp \left( \ Q^{**}(v/(1 - \epsilon)) \ \right) \ \right\}.
\end{equation}

\vspace{3mm}

 \ As a consequence:\par

\vspace{4mm}

 \ {\bf Proposition 2.3.} If

\begin{equation} \label{Y finite}
\exists \ \epsilon \in (0,1) \ \Rightarrow  Y_Q(\epsilon) < \infty,
\end{equation}
then the estimations (\ref{saddle}) holds true.  If in addition the function $ \ Q = Q(v) \ $ is continuous and convex,
then   (\ref{prelim Mr lambda}) is valid, as well. \par

 \vspace{5mm}

 \begin{center}

  {\sc To summarize.}

 \end{center}

  \hspace{3mm} {\bf A.} Recall that

 $$
 \ln M_f \left( \ e^v \ \right) \le \Lambda(v), \ v \in R,
 $$
 therefore

\begin{equation} \label{cn estim}
|c_n[f]| \le \exp \left( \ - \Lambda^*(n)  \  \right), \ n = 1,2,\ldots.
\end{equation}

\vspace{4mm}

\ {\bf B.}  Conversely, let the inequality (\ref{cn estim}) be given for certain non - negative continuous convex function
$ \ \Lambda = \Lambda(v), \ v \in R \ $  for which

$$
\exists \epsilon_0 \in (0,1) \ \Rightarrow S_0 = S(\epsilon_0) \stackrel{def}{=} Y_{\Lambda^*(\epsilon_0)} < \infty.
$$

 \ Then

\begin{equation} \label{Mr estim}
 M_f(r) \le S(\epsilon_0) \ \exp \left\{ \Lambda \left( \ \frac{  \ln r}{ 1 - \epsilon_0} \ \right) \ \right\}, \ r \ge e.
\end{equation}

\vspace{3mm}

 \ Briefly: under formulated above conditions

\begin{equation} \label{brief 1}
M_f(r) \le e^{ \ \Lambda(\ln r)  \ } \ \Rightarrow |c_n| \le e^{-\Lambda^*(n)}, \ r \ge e;
\end{equation}

\begin{equation} \label{brief 2}
|c_n| \le e^{-\Lambda^*(n)} \ \Rightarrow M_f(r) \le S_0 \ e^{ \ \Lambda(\ln r/( \ 1 - \epsilon_0 \ )) \ }.
\end{equation}

\vspace{4mm}

\ {\bf C. \ "Tauberian" theorem. } \par

 \ Let us impose an {\it additional} restriction on the function $ \ \Lambda(\cdot): \ $

\begin{equation} \label{gamma cond}
\exists \gamma =\gamma(\Lambda,\epsilon_0) = \const < \infty \ \Rightarrow \Lambda \left( \ \frac{v}{1 - \epsilon_0} \ \right) \le \gamma \ \Lambda(v), \ v \ge 1.
\end{equation}

 \ It follows immediately from the relations (\ref{brief 1}) and (\ref{brief 2}) the following assertion. \par

\vspace{4mm}

{\bf Theorem 2.1.} We conclude under formulated above conditions

\begin{equation} \label{Tauberian}
\lim_{r \to \infty} \frac{\ln M_f(r)}{\Lambda(\ln r) } = \lim_{n \to \infty} \frac{|\ln 1/|c_n||}{\Lambda^*(n)}.
\end{equation}

 \ More precisely: if there exists the left - hand side of (\ref{Tauberian}), then there exists also the right - hand one and they are equal;
  the converse proposition is also true: if there exists the  right - hand side of (\ref{Tauberian}), then there exists also the left - hand one and they are equal.\par

\vspace{5mm}

 \section{Examples}

\vspace{4mm}

 \ {\it Auxiliary fact.} Introduce a following family of regular varying functions

 \begin{equation} \label{family}
 \phi_{m,L}(\lambda) \stackrel{def}{=} \frac{1}{m} \ \lambda^m \ L(\lambda), \ \lambda \ge 1.
 \end{equation}
 \ Here $ \ m = \const > 1, \ $  and define as ordinary $ \ m' := m/(m-1) \ $ and $ \ L = L(\lambda), \ \lambda \ge 1 \ $
 is positive continuous {\it slowly varying } as $ \ \lambda \to \infty \ $ function. It is known that as $ \  x \to \infty, \ x \ge 1 \ $

\begin{equation} \label{customary}
\phi^*_{m,L} \sim (m')^{-1} \ x^{m'} \ L^{ \ - 1/(m-1)   \ } \left( \ x^{1/(m-1)} \ \right),
\end{equation}
see \cite{Seneta}, pp. 40 - 44; \cite{Kozachenko}. For instance, if $ \ \phi_m(\lambda) = m^{-1} |\lambda|^m, \ \lambda \in R, \ $ then

$$
\phi_m^{*}(x) = (m')^{-1} \ |x|^{m'}, \ x \in R.
$$

\vspace{4mm}

 \ More generally, if

 $$
   \psi_{m,q}(\lambda) = C_1 \lambda^m \ [\ln \lambda]^q, \ \lambda \ge e, \ m = \const > 1, \ q = \const \ge 0, \
 $$
then

$$
\psi^*_{m,q}(x) \asymp C_2(m,q) \  x^{m'} \ [\ln x]^{-q/(m-1)}, \ x \ge e.
$$

\vspace{4mm}

 \ {\bf  Example 3.1.} Suppose that for some entire function $ \ f = f(z) \ $

\begin{equation} \label{m case}
\ln M_f(r) \asymp C_3(m) \  [\ln r]^m, \ r \ge e, \ m = \const > 1.
\end{equation}

 \ Then

\begin{equation} \label{coeff m}
c_n[f] \le \exp \left\{ \  - C_4(m) \ n^{m'}   \ \right\}, \ n \ge 0;
\end{equation}
and conversely  proposition is also true: from the estimation (\ref{coeff m}) follows the inequality (\ref{m case}). \par

 \vspace{3mm}

 \ Note that the case $ \ m \le 1 \ $ is trivial: the function $ \ f(z) \ $ is polynomial, see
 \cite{Levin}, chapter 1, sections 2 - 4. This implies that

$$
\exists N \in \{ \ 2,3,\ldots \ \} \hspace{3mm}  \forall n \ge N \ \Rightarrow c_n = 0.
$$

\vspace{4mm}

\ {\bf  Example 3.2.} Suppose that for some entire function $ \ f = f(z) \ $

\vspace{3mm}

\begin{equation} \label{rho case}
\ln M_f(r) \sim C_4 \  r^{\rho}, \ r \ge 1, \ \rho = \const > 0, \ C_4 = \const \in (0,\infty),
\end{equation}
a classical case, \  \cite{Levin}, chapter 1, sections 1 - 5. \ Then

\begin{equation} \label{coeff rho}
|c_n| \le   \left[ \ \frac{n}{C_4 \ \rho} \ \right]^{-n/\rho} \cdot e^{n/\rho},
\end{equation}
 and conversely  proposition is also  true: from the estimation (\ref{coeff rho}) follows the inequality (\ref{rho case}). \par

  \ More  generally, the relation of the form

$$
\ln M_f(r)  \sim \rho^{-1} r^{\rho} \ \ln^{\gamma}r, \ r \to \infty
$$
is quite equivalent to the following  equality:  as $ \ n \to \infty \ $

$$
\ln \left( \ \frac{1}{|c_n|}  \ \right) \sim \rho^{-1} n \ln n + \gamma n \ln \ln n/\rho  - \frac{n}{\rho}.
$$

\vspace{4mm}

\ {\bf  Example 3.3.} We propose that for arbitrary entire function $ \ f = f(z) \ $ the following relations are
equivalent:

\begin{equation} \label{exp level}
\exists C_5, C_6 \in (0,\infty) \ \Rightarrow \ln M_f(r) \le C_5 e^{ \ C_6 r \ }, \ r \ge 0,
\end{equation}

\vspace{3mm}

and

\vspace{3mm}

\begin{equation} \label{log cn}
\exists C_7 \in (0,\infty) \ \Rightarrow  |c_n| \le C_7 (\ln n)^{-n}, \ n \ge 3.
\end{equation}

\vspace{5mm}

 \section{Generalization on the entire (holomorphic)  functions of several complex variables.}

\vspace{4mm}

 \hspace{3mm}  It is no hard  to generalize  the obtained results on the case of the  analytical functions $ \ f = f(z) \ $ of
 several complex variables. \par
 \ Let us introduce first all some (ordinary) used notations. The (finite) dimension of a considered problem  will be denoted
 by $ \ d; \ d = 2,3,\ldots. \ $
 Correspondingly the multivariate variable $ \ z \ $ consists on the $ \ d \ $ independent complex variables

$$
z = \vec{z} = \{ \ z_1,z_2, \ldots, z_d \ \}.
$$

 \ Ordinary vector notations:

$$
k = \vec{k} = \{k_1,k_2, \ldots,k_d \}, \ k_j = 0,1,2,\ldots;  \ |k| := \sum_{j=1}^d k_j;
$$

$$
\vec{z}^{\vec{k}} = z^k = \prod_{j=1}^d z_j^{k_j};
$$

$$
\vec{r} = r  = \{r_1,r_2, \ldots, r_d\} \in R_+^d, \ r_j \ge 0;
$$

$$
v = \vec{v} = \vec{v}(r) := \exp(\vec{r}) = \{\ e^{r_1}, e^{r_2}, \ldots, e^{r_d}  \ \} \in R^d,
$$

$$
\Longleftrightarrow  \ln \vec{v} = \ln v   = \{ \ \ln v_j \ \}, \ j = 1,2,\ldots,d.
$$

\vspace{3mm}

 \ The {\it multivariate} Young - Fenchel  transform $ \ g^*(y), \ y \in R^d \ $  for the
 given numerical valued function $ \ g = g(x), \ x \in R^d \ $ having the  convex non - empty domain of definition $ \ \Dom[g] \ $ is defined
 as usually

\begin{equation} \label{Young Fen multi}
g^*(y) \stackrel{def}{=} \sup_{x \in \Dom[g]} (  \  (x,y) - g(x) \ ),
\end{equation}
where $ \  (x,y)   \ $ denotes the inner, or scalar product of the two $ \ d - \ $ dimensional vectors $ \ x, y. \ $ As above, if the function
$ \ g = g(x) \ $ is continuous and convex, $ \ g^{**} = g. \ $ The famous Young inequality has a form

\begin{equation} \label{Young in mult}
(x,y) \le g(\gamma x) + g^*(y/\gamma), \  \gamma = \const > 0, \ x,y \in R^d.
\end{equation}

\vspace{3mm}

 \ Further, the analytical  function $ \ f = f(z) \ $  has a form

\begin{equation} \label{multi f}
f(\vec{z}) = f(z) = \sum_{k} c_k \ z^k = \sum_{\vec{k}} c_{\vec{k}} \ \vec{z}^{\vec{k}},
\end{equation}
where  the numbers $ \ \{ c_k \} = \{ \ c_{\vec{k}} \ \}  \ $ are  the Taylor's coefficients for $ \ f = f(z) \ $
and the series in (\ref{multi f}) converges for all the complex vectors $ \ \vec{z}. \ $\par

 \ The maximal function $ \ M(r) = M_f(r) = M_f(\vec{r}) \ $ for $ \ f = f(z)  \ $ at the point $ \ r \in R^d_+ \ $ is defined as before

\begin{equation} \label{vec M}
M_f(\vec{r}) \stackrel{def}{=} \max_{|z_j| \le r_j} |f(\vec{z})| = \max_{|z_j| = r_j} |f(z)|.
\end{equation}

 \ Define also the function

\begin{equation} \label{Lamfun mult}
\Lambda(v) = \Lambda(\vec{v})  \stackrel{def}{=} \ln M \left( \ e^{\vec{v}}  \ \right) = \ln M_f(e^v), \ v \in R^d.
\end{equation}

\vspace{4mm}

\begin{center}

 \ {\sc  Upper estimate.} \par

\end{center}

\vspace{4mm}

\ We apply the following estimate

$$
|c_{\vec{k}}| \le \frac{M_f(\vec{r})}{{\vec{r}}^{\vec{k}}}, \  \vec{r} \in R^d_+,
$$
see, e.g., \cite{Korevaar}, page 15, formula 1.4.3. Following,

$$
|c_k| \le \exp \left( \  - ((k,v) - \Lambda(v)) \    \ \right);
$$

\vspace{3mm}

$$
|c_{\vec{k}}| \le \exp ( \  - \sup_{v \in R} (  (k,v) - \Lambda(v))  \  ),
$$

and as before \par

\vspace{4mm}

{\bf Proposition 4.1.}

\vspace{4mm}

\begin{equation} \label{lower mult}
|c_{\vec{k}}| \le \exp \left( \ - \Lambda^*(\vec{k})   \ \right).
\end{equation}

\vspace{4mm}

\begin{center}

 \ {\sc  Lower estimate.} \par

\end{center}

\vspace{4mm}

 \ The lower estimate is quite alike ones obtained in the second section in the one - dimensional case, i.e. when  $ \ d = 1; \ $
 we will apply also the methods explained in  \cite{Kozachenko}. \ Suppose

\begin{equation} \label{lower ck mult}
|c_k| =   |c_{\vec{k}}| \le \exp \left( \ - Q(\vec{k})  \  \right) = \exp \left( \ - Q(k) \ \right)
\end{equation}
for  certain increasing to infinity relative the all variables $ \ k_j, \ j = 1,2,\ldots,d \ $
function $ \ Q = Q(z), \  \ z = \vec{z} = \{z_j \}, \ z_j  \ge 0. \ $   We have for all the sufficient greatest values $ \ r_j \ge e \ $

\begin{equation} \label{Mr est1 mult}
M_f(\vec{r}) \le  \sum_{ \vec{k} \ge \vec{0} } \exp \left( \ (k,v) - Q(k)  \ \right) = : R(v) = R(\vec{v}), \ \vec{v} =  \vec{\ln r}.
\end{equation}

 \vspace{4mm}

 \  Define the following function

$$
K(\epsilon) := \sum_{\vec{k} \ge \vec {0}}   \exp \left( \ -\epsilon Q^*(\vec{k}) \ \right), \ \epsilon \in (0,1).
$$

 \ It is proved in  particular in  \cite{Kozachenko} that if $ \ \exists \epsilon \in (0,1) \ \Rightarrow K(\epsilon) < \infty,  \  $ then

$$
R(v) \le K(\epsilon) \exp \left( \ (1 - \epsilon) Q^{**} \left(  \  \frac{v}{1 - \epsilon} \ \right)    \ \right).
$$
 \ Further, as long as

$$
(k,v) \le Q^*( ( 1 - \epsilon) \ v) + Q^{**}(k/(1 - \epsilon)), \ \epsilon \in (0,1),
$$
we have

$$
R(v) \le U(\epsilon) \exp \left( \ Q^{**}(v/(1 - \epsilon)) \ \right),
$$
where

$$
U(\epsilon) := \sum_{\vec{k} \ge \vec{0}} \exp (Q^*((1 - \epsilon) \ \vec{k}) - Q^*(\vec{k})).
$$

\vspace{3mm}

\ To summarize:

\vspace{4mm}

 \ {\bf Proposition 4.2.} Denote

$$
Y(\epsilon) := \min(U(\epsilon), K(\epsilon)), \ \epsilon \in (0,1).
$$

 \ We assert

\begin{equation} \label{General lower}
R(v) \le Y(\epsilon) \ \exp \left( \ Q^{**}(v/(1 - \epsilon)) \ \right),  \ \epsilon \in (0,1).
\end{equation}

 \ Of course,

\begin{equation} \label{General lower  inf}
R(v) \le  \inf_{\epsilon \in (0,1)} \left\{ \ Y(\epsilon) \ \exp \left( \ Q^{**}(v/(1 - \epsilon)) \ \right) \ \right\}.
\end{equation}

\vspace{3mm}

 \ As a consequence:\par

\vspace{4mm}

 \ {\bf Proposition  4.3.} If

\begin{equation} \label{Y finite}
\exists \ \epsilon \in (0,1) \ \Rightarrow  Y(\epsilon) < \infty,
\end{equation}
then the following estimation

\begin{equation} \label{saddle  mult}
M_f(\vec{r}) \le \exp( \ Q^{**}(C \ \ln r) \ ), \  r_j \ge e.
\end{equation}
holds true.  \par

\vspace{3mm}

 \ If in addition the function $ \ Q = Q(\vec{v}) \ $ is continuous and convex,
then

\begin{equation} \label{saddle Q mult}
M_f(  \vec{r}) \le \exp( \ Q(C \ \ln r) \ ), \  r_j \ge e.
\end{equation}

 \ In particular, if the function $ \ Q(\cdot)  \ $  coincides with the introduced before function $ \ \Lambda^*(\cdot), \ $
  we conclude in this case

\begin{equation} \label{ultimate Mr}
M_f(r) \le  \exp \left( \  \Lambda(C \ \ln r)  \ \right),  \  r_j \ge e.
\end{equation}

\vspace{4mm}

\begin{center}

 \ {\sc  Multivariate examples.} \par

\end{center}

\vspace{4mm}

 \ It is no  difficult  to show the exactness of obtained estimates   still in the multidimensional case.  It is sufficient to
 consider the case $ \ d = 2 \ $ and the so - called factorizable function

$$
f = f(z) = f(z_1,z_2) = f_1(z_1) \cdot f_2(z_2), \ z = (z_1,z_2),
$$
where the functions $ \ f_1,f_2 \ $ are  function - examples considered  in the third section.  If

$$
f_1(z_1) = \sum_{k = 0}^{\infty} a_k \ z_1^k, \hspace{4mm} f_2(z_2) = \sum_{l=0}^{\infty} b_l \ z_2^l,
$$
then

$$
f(z_1,z_2) = \sum \ \sum_{k,l = 0}^{\infty} a_k  \ b_l \ z_1^k \ z_2^l= \sum \sum_{k,l = 0}^{\infty} \ c_{k,l} \ z_1^k \ z_2^l,
$$
where $ \ c_{k,l} = a_k \ b_l.  \ $ \par

  \ Obviously,

$$
M_f(r_1,r_2) = M_{f_1}(r_1) \cdot M_{f_2}(r_2) = \exp \left\{ \ \Lambda_{f_1}(\ln r_1) + \Lambda_{f_2}(\ln r_2) \ \right\}.
$$

 \ As long as

$$
\sup_{\lambda,\mu \in R} [ \ (\lambda x + \mu y)  - (g(\lambda) + h(\mu))  \ ] = g^*(x) + h^*(y),
$$
we conclude

$$
|c_{k,l}| \le \exp \left\{  \ - \Lambda^*_{f_1}(k) - \Lambda^*_{f_2}(l)  \   \right\}
$$
and so one; see theorem 3.1. \par

 \vspace{5mm}

 \section{Concluding remarks.}

 \vspace{5mm}

 \hspace{4mm} {\bf A.} It is no hard to generalize obtained estimations on the {\it derivatives } of the
 source function $ \ f = f(z), \ $ including the partial derivatives in the case of the function of several
 complex variables. For instance, for the function from (\ref{source function}) we have

$$
f'(z) = \sum_{k=1}^{\infty} \ k c_k[f] \ z^{k-1},
$$
so that

$$
c_k[f'] = (k+1) \ c_{k+1}[f].
$$

\vspace{3mm}

 \ {\bf B.} The continuous version of our estimations, i.e. the Tauberian theorems for the Laplace transform
are well known, see  e.g. \cite{Widder}. The non - asymptotical estimates may be found, e.g. in \cite{Ostrovsky0}, pp. 27 - 37. \par

\vspace{3mm}

 \ {\bf C.}  Let us list  briefly a several works devoted to the applications of Tauberian Theorems in the Probability Theory:
\cite{Bagdasarov}, \cite{Bingham}, \cite{Kruglov Yakimiv}, \cite{Nakagava1},  \cite{Nakagava2},  \cite{Nakagava3}.\par

\ Let us show some extension of obtained in these works results based on the our estimations. Let $ \ \xi \ $ be an integer
values non - negative random variable (r.v.):

\begin{equation} \label{rand}
{\bf P} (\xi = k) = c_k, \ k = 0,1,2,\ldots .
\end{equation}
 \ Of course, $ \ c_k \ge 0, \ \sum_k c_k = 1. \ $  The so - called {\it generating function} \ (g.f.) \ $ \  g[\xi] = g[\xi](z) \ $  for this r.v. is
 as ordinary defined as follows

\begin{equation} \label{generatrix}
g[\xi](z) \stackrel{def}{=}  {\bf E} z^{\xi} = \sum_{k=0}^{\infty} c_k \ z^k.
\end{equation}
 \ This notion play a very important role in the probability theory, in particular, in the reliability theory, in the grand deviation theory,
 in the theory of queue theory etc. It is important especially for these applications the asymptotical behavior of
 $ \ g[\xi](z) \ $ as $ \ |z| \to \infty. \ $  \par
 \ One can for example apply our theorem 2.1 for the probability theory. Namely, we conclude under formulated in this theorem conditions

\begin{equation} \label{Tauber probab}
\lim_{r \to \infty} \frac{\ln M_g(r)}{\Lambda_P(\ln r) } = \lim_{n \to \infty} \frac{|\ln 1/|c_n||}{\Lambda_P^*(n)}.
\end{equation}

 \ More precisely: if there exists the left - hand side of (\ref{Tauber probab}), then there exists also the right - hand one and they are equal;
  the converse proposition is also true: if there exists the  right - hand side of (\ref{Tauber probab}), then there exists also the left - hand
  one and they are equal.\par

 \ Here with accordance  (\ref{Lamfun})

\begin{equation} \label{Lamfun probab}
\Lambda_P(v) \stackrel{def}{=} \ln M_g(e^v), \ v \in R.
\end{equation}

 \ Analogous fact holds true still in the multidimensional case. \par

\vspace{6mm}

\vspace{0.5cm} \emph{Acknowledgement.} {\footnotesize The first
author has been partially supported by the Gruppo Nazionale per
l'Analisi Matematica, la Probabilit\`a e le loro Applicazioni
(GNAMPA) of the Istituto Nazionale di Alta Matematica (INdAM) and by
Universit\`a degli Studi di Napoli Parthenope through the project
\lq\lq sostegno alla Ricerca individuale\rq\rq (triennio 2015 - 2017)}.\par


\vspace{6mm}

\begin{thebibliography}{36}



\bibitem{Bagdasarov}
{\bf D.R.Bagdasarov, E.I.Ostrovsky.}  {\it Inversion of Chebyshev's inequality.}
Theory Probab. Appl., V. 40, issue 4, 873 - 878.



\bibitem{Bingham}
{\bf  H.Bingham.} {\it Tauberian theorems and large deviations.}
Stochastics,  {\bf 80}, (2008), 143 \ - \ 149 .

\bibitem{Evgrafov}
{\bf M.A.Evgrafov.} {\it Asymptotic Estimates and Entire Functions.}
Dover Books on Mathematics, Paperback – April 15, 2020.


\bibitem{Fedoryuk}
{\bf M.V.Fedoryuk.}  {\it The saddle-point method.}  Moscow, Science (Nauka), (1977).


\bibitem{Freund}
{\bf  M.Freund and E.G\"orlich.} {\it On the Relation between Maximum Modulus, Maximum Term, and Taylor Coefficients of an Entire Function.}
Journal of approximation Theory,
{\bf 43,}  194 \ - \ 203, (1985).



\bibitem{Korevaar}
{\bf Jaap Korevaar, Jan Wiegerinckj.} \ {\it Several Complex Variables.} j.o.o.wiegerinck@uva.nlversion of November {\bf 18,} 2011.
Korteweg-de Vries Institute for Mathematics Faculty of Science University of Amsterdam.



\bibitem{Kozachenko}
{\bf Kozachenko Yu.V., Ostrovsky E.} {\it Equivalencs between tails,Grand LebSpace and Orlicz norms  for random variables
without  Kramer's  condition. } Bulletin of Taras Shevchenko National University of Kyiv, series Physic and Mathematic,
2018, {\bf 4,},  220  - 231.



\bibitem{KrasRut}
{\bf Krasnoselsky M.A., Routisky Ya. B.} {\it Convex Functions and Orlicz Spaces.}
 P. Noordhoff Ltd, (1961), Groningen.


\bibitem{Kruglov Yakimiv}
{\bf V.M.Kruglov,   Yakimiv A.L.}  {\bf Probabilistic Applications of Tauberian Theorems.}
Teor. Veroyatnost. i Primenen., 2006, Volume 51, Issue 4, 822 \ - \ 824



\bibitem{Levin}
{\bf B.Ja.Levin.} {\it Distribution of zeros of entire functions. } Revised Edition. Americal Mathematical
Society, Rhode Island,
Translations of Mathematical Monographs, Translation from Russian, 1980, Volume 5.


\bibitem{Nakagava1}
 {\bf Nakagawa, K.}  {\it On the Exponential Decay Rate of the Tail of a Discrete Probability Distribution.}
  Stochastic Models, vol. 20, no.1, pp.31 \ - \ 42, 2004.

\bibitem{Nakagava2}
 {\bf Nakagawa, K.}   {\it Tail probability of random variable and Laplace transform.}
  Applicable Analysis, vol.84, no.5, pp. 499 \ - \ 522, May 2005.


\bibitem{Nakagava3}
{\bf Nakagawa, K.} {\it On the Series Expansion for the Stationary Probabilities of an M/D/1 Queue.}
Journal of the Operations Research Society of Japan, vol.48, no.2, pp. 111 \ - \ 122.




\bibitem{Ostrovsky0}
{\bf  Ostrovsky E.I. } (1999). {\it Exponential estimations for Random Fields and its
applications,} (in Russian). Moscow-Obninsk, OINPE. \\



\bibitem{RaoRen1}
{\bf Rao M.M., Ren Z.D.} {\it Theory of Orlicz Spaces.} Marcel Dekker Inc., 1991.
New York, Basel, Hong Kong.

\bibitem{RaoRen2}
 {\bf Rao M.M., Ren Z.D.} {\it Applications of Orlicz Spaces.} Marcel Dekker Inc., 2002.
New York, Basel, Hong Kong.


\bibitem{Rockafellar}
{\bf R. T. Rockafellar.} {\it Convex analysis.}  Princeton Mathematical Series, No. 28. Princeton University Press,
Princeton, N.J., 1970


\bibitem{Seneta}
{\bf Eugene Seneta.}  {\it Regularly Varying Functions.} Lectures Notes in Mathematics, 508, (1976).


\bibitem{Tiel}
{\bf J. van Tiel. } {\it Convex Analysis: An Introductory Text.}  Wiley, 1984

bibitem{Ostrovsky0}
{\bf  Ostrovsky E.I. } (1999). {\it Exponential estimations for Random Fields and its
applications,} (in Russian). Moscow-Obninsk, OINPE. \\


\bibitem{Widder}
{\bf Widder, David Vernon.}  (1941.) {\it The Laplace Transform.}
Princeton Mathematical Series, v. 6, Princeton University Press, MR 0005923



\end{thebibliography}
\end{document}